\documentclass[reqno, 11pt]{amsart}
\usepackage{graphicx} % Required for inserting images
\usepackage[colorlinks=true]{hyperref}
\hypersetup{urlcolor=blue, citecolor=blue,linkcolor=blue}
\usepackage{amsmath, amssymb, amsthm, latexsym}
\usepackage{breqn}

\newtheorem{theorem}{Theorem}[section]
\newtheorem{corollary}[theorem]{Corollary}
\newtheorem{lemma}[theorem]{Lemma}
\newtheorem{proposition}[theorem]{Proposition}

\newtheorem{remark}[theorem]{Remark}
\newtheorem{example}[theorem]{Example}

\numberwithin{equation}{section}

\newcommand{\R}{\mathbb R}	% set of reals
\newcommand{\U}{\mathrm{U}}

\newcommand{\diag}{\operatorname{diag}}

\newcommand{\ol}[1]{\overline{#1}}
\newcommand{\blue}[1]{\textcolor{blue}{#1}}

\def\Cnn{{\mathbb C}_{n\times n}}
\def\H{\mathbb H}
\def\P{\mathbb P}
\def\la{\lambda}

\begin{document}

\title{Order Relations of the Wasserstein mean and the spectral geometric mean}
\date{\today}

\author{Luyining Gan}
\address{School of Science, Beijing University of Posts and Telecommunications, Beijing 100876, China}
\email{elainegan@bupt.edu.cn}

\author{Huajun Huang}
\address{Department of
Mathematics and Statistics, Auburn University, Auburn, AL, 36849, USA}
\email{huanghu@auburn.edu}



\begin{abstract}
On the space of positive definite matrices, several operator means are popular and have been studied extensively. %, especially the Wasserstein mean.
In this paper, we investigate the near order and the L\"owner order relations on the curves defined by the Wasserstein mean and the spectral geometric mean. 
We show that the near order $\preceq $ is stronger than the eigenvalue entrywise order, and that
$A\natural_t B \preceq A\diamond_t B$ for $t\in [0,1]$. We prove the monotonicity properties of 
the curves originated from the Wasserstein mean and the spectral geometric mean in terms of the near order.  
 The L\"owner order properties of the Wasserstein mean and the spectral geometric mean are also explored.
\end{abstract}

\maketitle

\section{Introduction}

The primary objective  of this manuscript is to investigate the near order and the L\"owner order relations concerning  the Wasserstein mean and the spectral geometric mean, as well as  the curves induced by these two means. 

Let $\Cnn$ be the space of all $n \times n$ complex matrices. In  $\Cnn$, let $\H_n$ (resp. $\P_n$, $\ol{\P}_n$) be the set  of $n\times n$ Hermitian (resp.  positive definite, positive semidefinite) matrices, 
% $\ol{\P}_n$ the set of $n\times n$ positive semidefinite matrices in $\Cnn$,
and $\U(n)$  the group of $n\times n$ unitary matrices.
Given $A \in \Cnn$, we denote $|A| = (A^*A)^{1/2}$.
For $A\in \P_n$, let $\la(A)$   denote the $n$-tuple of eigenvalues of $A$ with nonincreasing order, that is, $\la(A) = (\la_1(A), \la_2(A),\dots, \la_n(A))$ and $\la_1(A) \geq \la_2(A) \geq \cdots \geq \la_n(A)$.

On the space of positive definite matrices, several operator means are popular and have been studied extensively.
Given $A, B\in \P_n$, \emph{the arithmetic mean}, 
\emph{the Wasserstein mean}, \emph{the metric geometric mean}, and \emph{the spectral geometric mean} are defined for $t\in[0,1]$:
\begin{eqnarray}
    A\nabla_t B &=& (1-t)A+tB,  \label{eqn:am}\\
     A\diamond_t B &=&(1-t)^2 A+ t^2 B+t(1-t)[(AB)^{1/2}+(BA)^{1/2}], \label{eqn:wm}\\
 A \sharp _t B  &=& A^{1/2}(A^{-1/2}BA^{-1/2})^{t}A^{1/2}, \label{eqn:gm}\\
A \natural_t B  &=& (A^{-1}\sharp  B)^{t}A(A^{-1}\sharp  B)^{t}, \label{eqn:sm}    
\end{eqnarray}
where $A\sharp B = A\sharp _{1/2} B$. 
Many discovered relations  between these means are related to
their spectra. 

 The metric geometric mean was first introduced by Pusz and Woronowicz~\cite{PW75} for $t=1/2$ and  extended to $t\in [0,1]$ by Kubo and Ando~\cite{KA79}. The   mean was extended  to multiple variables  by Ando, Li and Mathias~\cite{ALM}. The properties of this mean were studied extensively by Lim~\cite{Li12}.  
 The spectral geometric mean was introduced by Fiedler and Pt\'ak~\cite{FP97}  for $t=1/2$ and extended to  $t\in [0,1]$   by Lee and Lim~\cite{LL07}. Properties of spectral geometric mean can be found in~\cite{AKL07, GK23, GT22, KL07, Kim21} and the references therein. 
 The Wasserstein mean is linked to the barycenter  in the Wasserstein space of Gaussian distributions, which is one of the popular topics in matrix analysis and probability theory~\cite{AC11, ABCC16}. 
Many interesting inequalities and properties of the Wasserstein mean has been given in~\cite{BJL19, GK23, HK19, HK20}. 
 Hwang and Kim in~\cite[Lemma 2.4]{HK22} gave another form of the Wasserstein mean~\eqref{eqn:wm} in terms of the metric geometric mean and arithmetic mean as
\begin{equation}\label{eqn:Wasserstein}
    A\diamond_t B=[I\nabla_t(A^{-1}\sharp  B)]A[I\nabla_t(A^{-1}\sharp  B)].
\end{equation}
This form is analogous to the form of the spectral geometric mean~\eqref{eqn:sm}. In consequence, many properties of these two means are similar.
%We will extend some properties from one to the other. 

The following interesting relations between positive definite matrices will be discussed.
\begin{enumerate}
    \item  We  use $A \geq 0$ to denote that $A$ is positive semidefinite.    Two matrices $A, B \in \P_n$ satisfy \emph{the  L\"owner order} $A\le B$ if $B-A\ge 0$.  

\item 
 Define \emph{the near order} on $\P_n$: $A\preceq B$ if $A^{-1}\sharp  B\ge I$. The relation is introduced by Dumitru and Franco~\cite{DF23}.

\item 
Given $A, B \in \mathbb{P}_n$, we define the \emph{eigenvalue entrywise order} relation: $A \leq_{\la} B$ if $\la_i(A) \leq \la_i(B)$ for $1\leq i\leq n$.
We say that $A=_{\la} B$ if $\lambda(A)=\lambda(B)$.

    % \item Let $x, y$ be two $n$-tuples of positive real numbers. Denote by $x^{\downarrow}, y^{\downarrow}$ the nonincreasing order of elements of $x, y$ respectively.
\item Given $A, B \in \mathbb{P}_n$, we write $A \prec_{w \log} B$ if $\lambda(A)$ is \emph{weakly log-majorized} by $\lambda(B)$, that is,
\begin{equation}\label{eqn:wlog}
\prod_{i=1}^k \lambda_i(A)\leq \prod_{i=1}^k \lambda_i(B), \quad k = 1, 2, \dots, m.
\end{equation}
In particular, we say that $\lambda(A)$ is \emph{log-majorized} by $\lambda(B)$, denoted by $A \prec_{\log} B$, if \eqref{eqn:wlog} is true for $k = 1, 2, \dots, n-1$ and equality holds for $k=m$.

% Given $A, B \in \mathbb{P}_n$,  we write $A \prec_{w \log} B$  if $\lambda(A) \prec_{w \log} \lambda(B)$,
% and $A \prec_{ \log} B$ if $\lambda(A) \prec_{ \log} \lambda(B)$.

\end{enumerate}

The above relations satisfy the transitive property except for  the near order.
%The near order does not satisfy the transitive property, 
The L\"owner order is the strongest condition among these relations.
The near order is   weaker than the L\"owner order but stronger than the eigenvalue entrywise order, which we will prove in Theorem \ref{thm: near order eig rel}.
 It is straightforward  that $A \leq_{\la} B$ implies $A \prec_{w \log} B$.  
So for $A, B\in \P_n$, we have the following relationships:
\begin{equation}\label{eqn:relationship}
    A\le B \Longrightarrow A\preceq B \Longrightarrow A\le_{\la} B \Longrightarrow A\prec_{w\log} B.
\end{equation}
In Theorem \ref{thm:near-main}, we show that  $A\natural_t B \preceq A\diamond_t B$ for $t\in (0,1)$, 
which completes a chain  of order relations among   means in \eqref{eqn:chain}. 
Consequently, there is the  eigenvalue entrywise relation
 $A\natural_t B \le_{\lambda} A\diamond_t B$ for $t\in (0,1)$.

The metric geometric mean satisfies that if $A\le B$, then the induced geodesic curve 
 $\{A\sharp _{t} B\mid t\ge 0\}$ is monotonically increasing in terms of the  L\"owner order $\le$  with respect to $t$. Similar properties hold for the curves  induced by the Euclidean mean $A\nabla_t B$ and the log-Euclidean mean $\exp({(1-t) \log A + t \log B})$,   but not for the   curves induced by the Wasserstein mean and the spectral geometric mean. 
We show in Theorems \ref{thm: A preceq B}  that if $A\preceq B,$
then $\{A\diamond_{t} B\mid t\ge 0\}$  and $\{A\natural_{t} B\mid t\in \R \}$ are monotonically increasing in terms of the near order $\preceq$ with respect to $t$. 
In Theorem \ref{thm:near_order_range}, the near order relations between $A\diamond_t B$ and $A\diamond_s B$, between $A\natural_t B$ and $A\natural_s B$, and between $A\natural_t B$ and  $A\diamond_s B$,  are compared, respectively. For certain real powers $p$, 
  the near order  on the curves  $\{A^p\diamond_t B^p\mid t\ge 0\}$ and $\{B^{-p}\diamond_t A^{-p}\mid t\ge 0\}$ are also discussed. The results disclose the existence of abundant near order relations and the corresponding eigenvalue entrywise order relations in the Wasserstein mean, the spectral geometric mean, and the curves induced by  two means.

In Section~\ref{Sec:properties}, we study  the L\"owner order properties of Wasserstein mean and spectral geometric mean, both of which are connected to
the metric geometric mean. In particular, Theorem \ref{thm:was_geom} shows that if $A\diamond_t B\le A\diamond_t C$ or $A\natural_t B\le A\natural_t C$ for one $t\in(0, 1]$,  then $A^{-1}\sharp_s  B \le A^{-1}\sharp_s  C$ for all $s\in [0,1/2]$. 

Some  preliminary  results of 
the metric geometric mean, the Wasserstein mean, the spectral geometric mean, and their relations  
are explored in Section~\ref{Sec:pre}.

\section{Preliminaries}\label{Sec:pre}

Extensive investigations have been done on the properties of matrix means. 
Here we list some basic properties of 
the metric geometric mean~\cite{Bh07, KA79, LL07}, the spectral geometric mean~\cite{LL07, Li12} and the Wasserstein mean~\cite{HK19, HK22, KL20}. They display  different characteristics of these three means. 

\begin{theorem}\label{thm:gm}
    Let $A, B, C, D\in \mathbb{P}_n$ and let $s, u, t\in [0, 1]$. The following are satisfied.
    \begin{enumerate}
\item $A \sharp _t B = B \sharp _{1-t} A$.
\item $(A \sharp _t B)^{-1} = A^{-1} \sharp _t B^{-1}$.
\item $(aA) \sharp _t (bB) = a^{1-t} b^{t} (A \sharp _t B)$ for any $a, b >0$.
    \item $\det (A \sharp _{t} B) = (\det A)^{1-t} (\det B)^{t}$.
        \item $M(A\sharp _t B)M^* = (MAM^*)\sharp _t (MBM^*)$ for non-singular $M$.
        \item $(A^{-1}\nabla_t B^{-1})^{-1} \leq A\sharp _t B \leq A\nabla_t B$, where $A\nabla_t B = (1-t)A+tB$.
        % \item $(\lambda A+(1-\lambda)B)\sharp _t ((\lambda C+(1-\lambda)D)) \geq \lambda(A\sharp _t C) +(1-\lambda)(B\sharp _tD)$;
        \item $(A\sharp _s B)\sharp _t(A\sharp _u B) = A\sharp _{(1-t)s+tu} B$.
        % \item $\lim_{p\to 0} (A^p\sharp _t B^p)^{1/p} = \exp((1-t)\log A+t\log B)$.
        \item If $A\le C$ and $B\le D$, then $A\sharp_t B \le C\sharp_t D.$
        \item  $A^{-1}\sharp  (BAB)=B.$ 
    \end{enumerate}
\end{theorem}

\begin{theorem}\label{thm:sgm}
    Let $A, B\in \mathbb{P}_n$ and let $s, u, t\in [0, 1]$. The following are satisfied.
    \begin{enumerate}
\item $A \natural_t B = B \natural_{1-t} A$.
\item $(A \natural_t B)^{-1} = A^{-1} \natural_t B^{-1}$.
\item $(aA) \natural_t (bB) = a^{1-t} b^{t} (A \natural_t B)$ for any $a, b >0$.
\item $(A \natural_s B) \natural_t (A \natural_u B) = A \natural_{(1-t)s+tu} B$.
\item $\det (A \natural_{t} B) = (\det A)^{1-t} (\det B)^{t}$.
% \item $\lim_{p\to 0} (A^p\natural_t B^p)^{1/p} = \exp((1-t)\log A+t\log B)$.
    \end{enumerate}
\end{theorem}

%The definition \eqref{eqn:sm} of $A\natural_t B$ can be extended to all $t\in \R$. It is easy to check that all properties of the spectral geometric mean in Theorem \ref{thm:sgm} still hold for $s,u,t\in\R$. 

\begin{theorem} \label{thm:wm}
    Let $A, B \in \mathbb{P}_n$ and let $s, t, u \in [0,1]$. The following are satisfied.
\begin{enumerate}
\item $A \diamond_t B = B \diamond_{1-t} A$.
\item $(A \diamond_t B)^{-1} = A^{-1} \diamond_t B^{-1}$ if and only if $A = B$.
\item $(aA)\diamond_t (aB) = a (A \diamond_t B)$ for any $a > 0$.
\item $(A \diamond_s B) \diamond_t (A \diamond_u B) = A \diamond_{(1-t)s+tu} B$.
\item $\det (A \diamond_{t} B) \geq (\det A)^{1-t} (\det B)^{t}$.
\item $A\diamond_t B \leq A\nabla_t B$.
% \item $\lim_{p\to 0} (A^p\diamond_t B^p)^{1/p} = \exp((1-t)\log A+t\log B)$.
\end{enumerate}
\end{theorem}

Besides the individual properties of matrix means,  many relations among these means have been discovered. 
%% For the convenience of the reader, we collected some known results and our new results to show the relations between the matrix means.
The near order relation $\preceq$ on $\P_n$ is recently introduced by Dumitru and Franco~\cite{DF23}: 
\begin{center}
    $A\preceq B$ if and only if $A^{-1}\sharp B \geq I$.
\end{center} 
It is natural that $A\preceq B$ if and only if $B^{-1}\preceq A^{-1}$. However,
 the relation ``$\preceq$'' does not satisfy transitive property~\cite{DF23}, that is,
$A\preceq B$ and $B\preceq C$ do not necessarily imply that $A\preceq C$. So  the relation ``$\preceq$'' is called  a near order. 

On one hand, the near order relation is   weaker than the L\"owner order, since $A\le B$ implies that $A^{-1}\sharp B \ge A^{-1}\sharp A =I$ so that $A\preceq B$. 
On the other hand, we show below that the near order is stronger  than the  eigenvalue entrywise relation. % the (weak) log-majorization relation.

 \begin{theorem}\label{thm: near order eig rel}
    Let $A, B\in \P_n$. If $A\preceq B$, then $A\le_{\lambda} B$.
\end{theorem}
 
\begin{proof}
    Assume $A\preceq B$, we have $B^{-1}\sharp  A\le I$. Then it is straightforward to have
  \begin{eqnarray*}
      A &=_{\la}&  B^{1/2} B^{-1/2} (B^{1/2}A B^{1/2})^{1/2} B^{-1} (B^{1/2}A B^{1/2})^{1/2} B^{-1/2} B^{1/2} \\
        &=& B^{1/2} (B^{-1}\sharp  A)^2 B^{1/2} \\
        &\le& B^{1/2}IB^{1/2} = B. 
  \end{eqnarray*}
Thus $A\le_{\lambda} B$.
\end{proof} 

From~\cite{ AH94, A90, GLT21}, we have the log-majorization relations between the metric geometric mean, the log-Euclidean mean, the Fidelity mean and the spectral geometric mean as
\begin{equation}\label{eqn:means}
  A \sharp _{t} B \prec_{\log} \exp({(1-t) \log A + t \log B}) \prec_{\log} B^{t/2}A^{1-t}B^{t/2} \prec_{\log} A \natural_{t} B.   
\end{equation}
It is well-known  (e.g. \cite{LL11, BJL19} ) that there are L\"owner orders:
%between the metric geometric mean and the arithmetic mean and between the Wasserstein mean and the arithmetic mean, that is, 
\begin{equation}\label{eqn:meanlowner}
    (A^{-1}\nabla_t B^{-1})^{-1}\leq A\sharp _t B\leq A\nabla_t B \quad \text{and} \quad A \diamond_{t} B \le  A\nabla_t B.
\end{equation}
Our new result in Theorem \ref{thm:near-main}   shows that $A \natural_{t} B \preceq  A \diamond_{t} B$ for all $t\in (0,1)$. It completes a chain of order relations among   means as
\begin{eqnarray}\label{eqn:chain}
   (A^{-1}\nabla_t B^{-1})^{-1}\leq A \sharp _{t} B &\prec_{\log}& \exp({(1-t) \log A + t \log B}) \nonumber\\
   &\prec_{\log}& B^{t/2}A^{1-t}B^{t/2} \prec_{\log} A \natural_{t} B \preceq  A\diamond_t B \leq A\nabla_t B.   
   \qquad 
\end{eqnarray}

%Since the L\"owner order, the eigenvalue entrywise order, and the (weak) log-majorization order are transitive and the relations of these orders are given in~\eqref{eqn:relationship}, 
The relation between any two of the above means can be derived from   \eqref{eqn:chain}  and the transitive properties. 
Here are some relations:
\begin{enumerate}

    \item  \eqref{eqn:chain} and \cite[Theorem 2]{DF23}  imply that
    \begin{equation}
        A \natural_{t} B \preceq  A\nabla_t B.
    \end{equation}
    
    \item From \eqref{eqn:chain}, we get  
    \begin{equation}
        A \sharp _{t} B  \prec_{w \log} A\diamond_t B.
    \end{equation}
    It is the known strongest relation. The following counterexample shows that
    $A \sharp _{t} B$ and $A\diamond_t B$ do not have eigenvalue entrywise relation, near order relation, or L\"owner  order relation. 
    Let \[A = \begin{bmatrix}
 39.1195&   42.1116\\
   42.1116&   61.1568   
 \end{bmatrix} \quad \text{and} \quad 
 B = \begin{bmatrix}
        26.3279&   13.3485\\
   13.3485&   12.2727
 \end{bmatrix}.
\]
 Then we have
 \[A\sharp_t B = \begin{bmatrix}
        32.2446&   29.2497\\
   29.2497&   39.8872
 \end{bmatrix}
 \quad \text{and}\quad 
 A\diamond_t B = \begin{bmatrix}
        35.6339&   33.9111\\
   33.9111 &  45.3815
 \end{bmatrix}.
 \]
The spectra of $A\sharp_t B$ is $\{   65.5641, 6.5677\}$ and that of $A\diamond_t B$ is $\{74.7672, 6.2481\}$. 
Thus $A \sharp _{t} B$ and $A\diamond_t B$ do not have eigenvalue entrywise relation.  \eqref{eqn:relationship} implies that
$A \sharp _{t} B$ and $A\diamond_t B$  do not have either near order relation or L\"owner  order relation.

    \item  \eqref{eqn:chain} may not always provide the strongest relation between means. 
For example, \eqref{eqn:chain} implies that $(A^{-1}\nabla_t B^{-1})^{-1} \prec_{w \log} A\diamond_t B$. However,
Theorem~\ref{thm:wm}(6) and Theorem \ref{thm:wm_inverse} show that
\begin{equation}
(A^{-1}\nabla_t B^{-1})^{-1} \le (A^{-1}\diamond_t B^{-1})^{-1}  \preceq   A\diamond_t B\leq A\nabla_t B,    
\end{equation}
	hence by  \cite[Theorem 2]{DF23},    a stronger relation exists:
\begin{equation}
(A^{-1}\nabla_t B^{-1})^{-1} \preceq  A\diamond_t B.    
\end{equation}

\end{enumerate}

A relation between the near order and the L\"owner  order   is given as follows, which strengthen a result in \cite[Hwang and Kim, 2022]{HK22}.

\begin{theorem}\label{thm:Hermitian-near}
    Let $A\in\P_n$. Let $P,Q\in \H_n$ be nonsingular and $PQ=QP$. Then $I\le P^{-1}Q$ if and only if
    $PAP\preceq QAQ.$ 
\end{theorem}

\begin{proof} Since $P$ and $Q$ commute,  $P^{-1}Q=QP^{-1}\in\H _n.$
By Theorem \ref{thm:gm} (9),
$I\le P^{-1}Q$  if and only if 
$$I\le P^{-1}Q=(PAP)^{-1}\sharp [(P^{-1}Q)(PAP)(P^{-1}Q)]  =(PAP)^{-1}\sharp (QAQ),$$
 if and only if  $PAP\preceq QAQ.$
\end{proof}

\begin{corollary}\label{thm:near}
    Let $A, P,Q\in\P_n$ and $PQ=QP$, then  $P\le Q$ if and only if
    $PAP\preceq QAQ.$ 
\end{corollary}

By \cite[Theorem II.1]{ZH23}, on the manifold $\ol{P_n}$ of $n\times n$ positive semidefinite matrices, the formula  \eqref{eqn:wm} of $A\diamond_t B$  can be extended to  all  $t\in\R$,  and 
\begin{equation}\label{BW geodesic form}
A\diamond_t B =B\diamond_{1-t} A  = |(1-t)A^{1/2}+t U^* B^{1/2}|^2 =|A^{1/2}\nabla_t (U^* B^{1/2})|^2
\end{equation} 
in which $U$ is a unitary matrix in the polar decomposition $ B^{1/2}A^{1/2}=U |B^{1/2}A^{1/2}|$, and $\nabla_t$ in \eqref{eqn:am} is extended to all $t\in\R.$  When $A\in\P_n,$
\eqref{BW geodesic form} is equivalent to that for $t\in\R$: 
\begin{equation}\label{BW geodesic form 2} 
A\diamond_t B =[I\nabla_t(A^{-1}\sharp B)] A [ I\nabla_t(A^{-1}\sharp B)].
\end{equation} 
Likewise, the formula \eqref{eqn:sm} of $A\natural_t B$ can   be extended to all $t\in\R$ such that
\begin{equation}\label{spectral geometric form}
A \natural_t B  = (A^{-1}\sharp  B)^{t}A(A^{-1}\sharp  B)^{t}.  
\end{equation}

\eqref{BW geodesic form 2} and \eqref{spectral geometric form} immediately imply the following result.

\begin{proposition}
Let $A, B\in\P_n$. Let $\mu_1,\ldots,\mu_n$ be the eigenvalues of $A^{-1}\sharp B$. Then for $t\in\R$, 
\begin{eqnarray*}
\det(A\diamond_t B) &=&  \det(A) \prod_{i=1}^n (1-t+t\mu_i)^2,
\\
\det(A\natural_t B) &=&  \det(A)^{1-t}\det(B)^{t}=\det(A) \prod_{i=1}^n \mu_i^{2t}.
\end{eqnarray*}
\end{proposition}

All properties of
 the spectral geometric mean in Theorem \ref{thm:sgm} still hold for $s,u,t\in\R$, and they can be proved by \eqref{spectral geometric form}. The extension of Theorem \ref{thm:sgm}(4) to $s,u,t\in\R$ is proved below. 

\begin{proposition}\label{lm:sgm-geodesic-property}
    Let $A,B\in \P_n$. Then for $s,u,t\in\R$,
$$(A \natural_s B) \natural_t (A \natural_u B) = A \natural_{(1-t)s+tu} B. $$
\end{proposition}

\begin{proof}
    \eqref{spectral geometric form} and Theorem \ref{thm:gm}(9) imply that for $A,B\in\P_n$  and  $s\in\R$, 
$$A^{-1}\sharp (A\natural_s B)=A^{-1}\sharp [(A^{-1}\sharp B)^s A (A^{-1}\sharp B)^s ]= (A^{-1}\sharp B)^s,$$
so that for $t,s\in\R,$
\begin{equation*}
A\natural_t(A\natural_s B)=[(A^{-1}\sharp B)^s]^t A [(A^{-1}\sharp B)^s]^t= A\natural_{ts} B.         
\end{equation*}
Therefore, 
$$(A\natural_s B)\natural_{t} B=B\natural_{1-t}(B\natural_{1-s} A)= B\natural_{(1-t)(1-s)}A=A\natural_{(1-t)s+t} B,$$
so that
\begin{equation*}
    (A \natural_s B) \natural_t (A \natural_u B) 
= [A\natural_{\frac{s}{u}} (A \natural_u B)] \natural_t  (A\natural_u B)
= A \natural_{(1-t)\frac{s}{u}+t}  (A \natural_u B)=A \natural_{(1-t)s+tu} B.  \qedhere
\end{equation*} 
\end{proof}

However, the counterpart of Proposition \ref{lm:sgm-geodesic-property} for $A\diamond_t B$ is not true, since  the Wasserstein  curve $\{A\diamond_t B \mid t\in\R\}$
given by  \eqref{eqn:wm}, \eqref{BW geodesic form}, or \eqref{BW geodesic form 2} 
consists of several geodesic curves joint at some  boundary points  of $\ol{P_n}.$
A difference between $\natural_t$ and $\diamond_t $
lies in the fact  that: in \eqref{spectral geometric form}, $(A^{-1}\sharp B)^t$  for $t\in\R$ is always positive definite, 
but in   \eqref{BW geodesic form 2}, $I\nabla_t (A^{-1}\sharp B)$ is not. 

\begin{example}
  If $A,B\in\P_n$ are commuting, then  \eqref{BW geodesic form} shows that
$$A\diamond_t B=(A^{1/2}\nabla_t B^{1/2})^2=|A^{1/2}\nabla_t B^{1/2}|^2.$$ 
Let  $A=\diag(4,1)$ and $B=\diag(1,4)$. Then for $t, s\in\R $,  
    \begin{eqnarray*}
        A\diamond_{ts} B &=& \diag((2-ts)^2, (1+ts)^2)
        \\
        A\diamond_{t} (A\diamond_{s} B) &=& \diag([2(1-t)+|2-s|t]^2, (1-t+|1+s|t)^2).
    \end{eqnarray*}
 When $t\ne 0, 1$ and $s>2$, the $(1,1)$ entries of $A\diamond_{ts} B$ and $ A\diamond_{t} (A\diamond_{s} B)$ are different.
\end{example}

When $A$ and $B$ has near order relation,  the following properties hold for Wasserstein  curve.

\begin{proposition}\label{lem:wm-geodesic-property}
Suppose that $A, B \in \P_n$. 
\begin{enumerate}
    \item If $A\preceq B$. 
Then for $t, s\in [0,\infty)$,
\begin{equation*}\label{diamond geodesic property1}
    A\diamond_{t} (A\diamond_{s} B)= A\diamond_{ts} B.
\end{equation*}  
    \item If  $A\succeq B$. 
Then for $t, s\in (-\infty, 1]$,
\begin{equation*}\label{diamond geodesic property2}
    (A\diamond_{s} B)\diamond_{t} B= A\diamond_{(1-t)s+t} B.
\end{equation*}   
\end{enumerate}
\end{proposition}

\begin{proof}
\begin{enumerate}
    \item $A\preceq B$ implies that $A^{-1}\sharp B\ge I.$
For $ s\ge 0$ we have  $I\nabla_s(A^{-1}\sharp B)\ge I$. By Theorem \ref{thm:gm}(9),
$A^{-1}\sharp  (A\diamond_{s} B)= I\nabla_s(A^{-1}\sharp B)$, and for $t\ge 0$,
\begin{eqnarray*}
    A^{-1}\sharp  [A\diamond_{t} (A\diamond_{s} B)]
    &=&
    I\nabla_{t}[A^{-1}\sharp  (A\diamond_{s} B) ]
  = I\nabla_{t}[I\nabla_{s} (A^{-1}\sharp  B)]
  \notag \\
    &=& I\nabla_{ts} (A^{-1}\sharp  B)
    =  A^{-1}\sharp  (A\diamond_{ts} B).
\end{eqnarray*}
Therefore, when $t,s\ge 0$,  we have
$A\diamond_{t} (A\diamond_{s} B)= A\diamond_{ts} B. $

\item It can be proved by using $A\diamond_t B=B\diamond_{1-t}A$ and the preceding result.
\qedhere 
\end{enumerate}
\end{proof}

In the coming sections, we will focus on exploring  the near order relation  and the L\"{o}wner order  relations  between the spectral geometric mean and the  Wasserstein mean, and these relations on the curves defined by two means.

\iffalse
\blue{The following lemmas may not be used.} 

\begin{lemma}\cite[Hwang and Kim, 2022]{HK22}\label{lem:HK22}
    Let $A, B\in \P_n$. %%If $A\diamond_t B\le A$ for some $t\in (0,1]$, then $A^{-p}\sharp  B^p \le I$ for $p\ge 1$.
    If $A\diamond_t B\le B$ for some $t\in [0,1)$, then $A^{-p}\sharp  B^p \ge I$ for $p\ge 1$.
\end{lemma}

\begin{lemma}\cite[Hwang and Kim, 2022]{HK22}\label{lem:sxs}
    Let $S, X \in \P_n$. Then $SXS\le X$ implies $S\le I$. Furthermore, $SXS=X$ for $S, X\in \P_n$ if and only if $S=I$.
\end{lemma}

\begin{lemma}\label{lem:AAB}
    Let $A, B\in \P_n$. If $A\le A\diamond_t B$ for some $t\in (0,1]$, then $A^{-p}\sharp  B^p \ge I$ for $p\ge 1$.
\end{lemma}

\begin{proof}
    Since $A\le A\diamond_t B$ for some $[0,1)$, it is equivalent to 
    \[[I\nabla_t(A^{-1}\sharp  B)]^{-1} A [I\nabla_t(A^{-1}\sharp  B)]^{-1} \le A.\]
    By Lemma~\ref{lem:sxs} we have $(I\nabla_t(A^{-1}\sharp  B))^{-1}\le I$ and thus $A^{-1}\sharp  B\ge I$. 
    Taking the inverse map, we get $A\sharp  B^{-1} = (A^{-1}\sharp  B)^{-1}\le I$.
    For all $p\ge 1$, by the Ando-Hiai inequality $A^p\sharp B^{-p} \le I$, that is, $A^{-p}\sharp B^{p} \ge I$.
\end{proof}
\fi

\section{The near orders on Wasserstein  curves and spectral geometric curves}\label{Sec:new_relation}

%% Recently, Dumitru and Franco~\cite{DF23} defined a new relation $\preceq$ on $\P_n$ by  \begin{center}    $A\preceq B$ if and only if $A^{-1}\sharp B \geq I$. \end{center} 
There are abundant near order relations  on the curves originated from the Wasserstein  mean and the spectral geometric mean. 
We will explore these relations here.  

 $A\diamond_t B$ and   $A\natural_t B$ have the following near order relations.

\begin{theorem}\label{thm:near-main}
Suppose $A, B \in \P_n$.
\begin{enumerate}
    \item For $t\in (0,1)$, the spectral geometric mean and the  Wasserstein mean satisfy that
    \begin{equation}\label{near-order-means}
        A\natural_t B \preceq A\diamond_t B.
    \end{equation}
    \item If $A\preceq B$, then for $t\in  (1,\infty)$,
    \begin{equation}\label{near-order-right}
        A\natural_t B \succeq A\diamond_t B. 
    \end{equation}

    \item If $A\succeq B$, then for $t\in  (-\infty,0)$,
    \begin{equation}\label{near-order-left}
        A\natural_t B \succeq A\diamond_t B. 
    \end{equation}
    
\end{enumerate}
Moreover, any  equality in  \eqref{near-order-means},  \eqref{near-order-right}, or \eqref{near-order-left}  holds  if and only if $A= B$.
\end{theorem}
 
\begin{proof}
Given $a>0$, define the function $f_a(t) =a^t-1+t-at$. If $a=1$ then $f_1(t)=0$ for all $t\in\R$. Otherwise, we have $f_a(0)=f_a(1)=0$ and
$f_a''(t)=(\ln a)^2 a^t>0$ for all $t\in\R$. Hence $f_a(t)< 0$ for $t\in (0,1)$ and $f_a(t)> 0$ for $t\in (-\infty, 0)\cup (1,\infty).$

Let $X=A^{-1}\sharp  B$. Let $\mu_1\ge \cdots\ge \mu_n$ be the eigenvalues of $X$.  Then $X> 0$ and  $X^t$ commutes with $I\nabla_t X$ for $t\in\R.$  
Therefore, $X^t$ and $I\nabla_t X$ are simultaneously unitarily diagonlizable and the eigenvalues of 
$X^t- I\nabla_t X$ are $f_{\mu_i}(t)$ for $i=1,2,\ldots,n.$
The preceding argument shows that: 
\begin{enumerate}
    \item For $t\in (0,1)$, we have the L\"{o}wner order $0< X^t\le I\nabla_t X$. By  Corollary~\ref{thm:near}, 
    $$A\natural_t B=X^t A X^t \preceq (I\nabla_t X)A(I\nabla_t X)= A\diamond_t B.$$
    \item Suppose that $A\preceq B$. Then all $\mu_i\ge 1$. So for $t\in  (1,\infty),$ we have $0< I\nabla_t X\le X^t.$
     By Corollary~\ref{thm:near}, we get $A\natural_t B\succeq A\diamond_t B.$
     \item Suppose that  $A\succeq B$. Then all $\mu_i\in (0,1]$. For $t\in (-\infty,0),$ we have $0< I\nabla_t X\le X^t.$
      Corollary~\ref{thm:near} implies that $A\natural_t B\succeq A\diamond_t B.$
\end{enumerate}    
The values of $f_{\mu_i}(t)$ show  that the equality in  \eqref{near-order-means}, \eqref{near-order-right}, or \eqref{near-order-left}  holds if and only if all $\mu_i=1$, that is, $A=B.$
\end{proof}

\begin{example}
For $t\in (0,1)$, we have the near order $A\natural_t B\preceq A\diamond_t B$. However, they don't satisfy the stronger L\"owner order. Here is an counterexample:
    \[A=\begin{bmatrix}
        50 & 0\\0 & 10
    \end{bmatrix} \quad \text{ and }\quad
    B = \begin{bmatrix}
        57.8906 & 19.8885\\19.8885 & 62.1094
    \end{bmatrix}\]
    when $t=1/2$, the spectrum of $A\diamond_t B -A \natural_t B$ is $\{-0.21, 6.3338\}$. 
\end{example}

Theorem \ref{thm: near order eig rel} shows that the near order is stronger  than the  eigenvalue entrywise relation. 
So \eqref{near-order-means} implies the following result.

\begin{corollary}\label{thm:main}
    If $A, B \in \P_n$, then for $t\in (0,1)$,
    \begin{equation}
        A\natural_t B \leq_{\la} A\diamond_t B,
    \end{equation}
    where the equality holds if and only if $A= B$.
\end{corollary}

Corollary \ref{thm:main} strengthens the weak log-majorization result  in~\cite{GK23}:  for $t\in (0,1)$, $A\natural_t B \prec_{w \log} A\diamond_t B$.

\begin{corollary}\label{cor:natural_pre}
    For $t\in (0,1),$ we have $A\natural_t B\preceq A\nabla_t B.$ 
\end{corollary}

\begin{proof}
It is known~\cite{LL11, BJL19} that  $A\diamond_t B\le A\nabla_t B$ for $t\in (0,1)$, and  \eqref{near-order-means} shows that $A\natural_t B\preceq A\diamond_t B$.   By ~\cite{DF23}, we get $A\natural_t B\preceq A\nabla_t B.$ 
\end{proof}

% Corollary~\ref{cor:natural_pre} also implies that $A\natural_t B\leq_{\la} A\nabla_t B$ for $t\in (0,1)$ by Theorem~\ref{thm: near order eig rel}.

Theorem \ref{thm:wm}(5) shows that for $t\in (0,1),$ we have $\det(A\diamond_t B)\ge (\det A)^{1-t} (\det B)^{t}=\det(A\natural_t B)$. 
The following is a direct consequence of Theorem \ref{thm:near-main}. 

\begin{corollary}
 If $A\preceq B$ and $t\in  (1,\infty)$, or
 $A\succeq B$ and $t\in  (-\infty,0)$, then   
 $$\det(A\diamond_t B)\le (\det A)^{1-t} (\det B)^{t}.$$
\end{corollary}

Next, we study the monotonicity of  near order on the curves defined by $\diamond$ and $\natural$. 
 
\begin{theorem}\label{thm: A preceq B} 
The following are equivalent for $A, B\in \P_n$ and $t,s\in\R $:
\begin{enumerate}
    \item[(1)]  $A\preceq B$;
    \item[(2)]  $A\diamond_t B\preceq A\diamond_s B$ for certain $0\le t<s\le 1$;
    \item[(3)]  $A\natural_t B\preceq A\natural_s B$ for certain $ t<s $.
\end{enumerate}
 Moreover, when $A\preceq B$, the parametric curves
$\{A\diamond_{t} B\mid t\ge 0\}$ and $\{ A\natural_t B\mid t\in\R \}$ are   monotonically increasing  with respect to   the near order, that is, 
\begin{eqnarray*}
0\le t< s &\Longrightarrow & A\diamond_t B\preceq A\diamond_s B,
\\
t<s &\Longrightarrow & A\natural_t B\preceq A\natural_s B.
\end{eqnarray*}
\end{theorem}

\begin{proof}
 Obviously (1) implies (2) and (3). 
 
Suppose that (2) holds, namely, there are $t, s\in [0,1]$ with $t<s$ such that  $A\diamond_t B\preceq A\diamond_s B$. 
By \eqref{eqn:Wasserstein}, 
$$
[I\nabla_t(A^{-1}\sharp  B)]A[I\nabla_t(A^{-1}\sharp  B)]\preceq [I\nabla_s(A^{-1}\sharp  B)]A[I\nabla_s(A^{-1}\sharp  B)]. 
$$
Note that $I\nabla_t(A^{-1}\sharp  B)$ and $I\nabla_s (A^{-1}\sharp  B)$ are in $\P_n$ and they commute. 
By Corollary \ref{thm:near}, 
$$I\nabla_t(A^{-1}\sharp  B)\le I\nabla_s(A^{-1}\sharp  B).$$
So  $I\le A^{-1}\sharp B$, and thus $A\preceq B.$ We get (1).

A similar argument shows that (3) implies (1). 

Now suppose that  $A\preceq B.$ Then $A^{-1}\sharp B\ge I$. For any $0\le t<s$, we have %the L\"{o}wner order
$$
0< I\le I\nabla_t(A^{-1}\sharp  B)\le I\nabla_s(A^{-1}\sharp  B),
$$
so that by \eqref{BW geodesic form 2} and Corollary \ref{thm:near},
$$
A\diamond_t B=[I\nabla_t(A^{-1}\sharp  B)]A[I\nabla_t(A^{-1}\sharp  B)]\preceq [I\nabla_s(A^{-1}\sharp  B)]A[I\nabla_s(A^{-1}\sharp  B)]
=A\diamond_s B. 
$$
Similary, for any $t<s$, we have $0\le (A^{-1}\sharp  B)^t\le (A^{-1}\sharp  B)^s$, so that by \eqref{spectral geometric form}  and Corollary \ref{thm:near},
\begin{equation*}
   A\natural_t B=(A^{-1}\sharp  B)^t A (A^{-1}\sharp  B)^t \preceq (A^{-1}\sharp  B)^s A (A^{-1}\sharp  B)^s =A\natural_s B. \qedhere  
\end{equation*}
\end{proof}

Likewise, Theorem \ref{thm: A preceq B} has a counterpart theorem for $A\succeq B$ and the proof can be obtained by using
$A\diamond_t B=B\diamond_{1-t} A$ and
$A\natural_t B=B\natural_{1-t} A.$
Indeed, more precise comparisons can be done by Theorem \ref{thm:Hermitian-near} together with \eqref{BW geodesic form 2} and \eqref{spectral geometric form} as follows. 

\begin{theorem}\label{thm:near_order_range}
Suppose that $A, B\in \P_n$ are distinct. Let   $\mu_1\ge\cdots\ge\mu_n$ be the eigenvalues of $A^{-1}\sharp B$. 
\begin{enumerate}
    \item   For real numbers $t<s$,
    \begin{equation}\label{natural:s-t-range}
        \mu_n^{2s-2t} (A\natural_t B)\preceq A\natural_s B \preceq \mu_1^{2s-2t} (A\natural_t B).
    \end{equation}
 
    \item  For real numbers $t<s$, 
    if one of the following cases occurs:
     \begin{enumerate}
        \item $t,s\in [0,1]$, or
        \item $A\preceq B$ or $A\succeq B$, $\mu_1\ne 1$, and  $t,s\in (\frac{1}{1-\mu_1},\infty)$,  or
        \item 
        $A\preceq B$ or $A\succeq B$, $\mu_n\ne 1$, and $t,s \in (-\infty,\frac{1}{1-\mu_n})$,
    \end{enumerate}    
    then 
   \begin{equation}\label{ws:s-t-range}
             \left(\frac{1-s+s\mu_n}{1-t+t\mu_n}\right)^2 (A\diamond_t B) \preceq A\diamond_s B\preceq \left(\frac{1-s+s\mu_1}{1-t+t\mu_1}\right)^2 (A\diamond_t B).
    \end{equation}
 
    \item %Let $g_{s,t}(\mu)=\mu^{-t}(1-s+s\mu)$. 
    For $t,s\in\R$ such that $1-s+s\mu_1>0$ and $1-s+s\mu_n>0$, let
    $$
    m_{s,t}=\min\left \{\frac{1-s+s\mu_i}{\mu_i^t}\mid i\in[n]\right \},
    \quad M_{s,t}=\max\left \{\frac{1-s+s\mu_i}{\mu_i^t} \mid i\in[n]\right \},
    $$ 
    then  
    \begin{equation}\label{sg-ws:s-t-range}
        m_{s,t}^2 (A\natural_t B)\preceq A\diamond_s B \preceq M_{s,t}^2 (A\natural_t B).
    \end{equation}
\end{enumerate}
\end{theorem}

\begin{proof}
    % We prove the first statement. The others are proved similarly.
\begin{enumerate}
    \item For $t<s$, 
    \begin{equation*}
    A\natural_s B = (A^{-1}\sharp B)^s A (A^{-1}\sharp B)^s = (A^{-1}\sharp B)^{s-t} (A\natural_t B) (A^{-1}\sharp B)^{s-t} .
\end{equation*}
Since $\mu_n^{s-t}I\le (A^{-1}\sharp B)^{s-t}\le \mu_1^{s-t}I$, Corollary \ref{thm:near} implies \eqref{natural:s-t-range}. 
    \item Given   $t<s$, we define the function  
    \begin{equation}\label{f_st}
            f_{s,t}(\mu)=\frac{1-s+s\mu}{1-t+t\mu} =\begin{cases}
        1-s+s\mu, &t=0 
        \\
        \frac{s}{t}-\frac{s-t}{t^2(\mu-1+1/t)}, &t\ne 0,\ \mu\ne 1-\frac{1}{t}
    \end{cases}
    \end{equation}
    Denote
    \begin{eqnarray*}
        l_{s,t}=\min\left \{f_{s,t}(\mu_i)\mid i\in[n]\right \},
        \qquad
        L_{s,t}=\max\left \{f_{s,t}(\mu_i)\mid i\in[n]\right \}.
    \end{eqnarray*}
    When $l_{s,t}>0$, we get
   $$
  0<l_{s,t}  I \le [I\nabla_s (A^{-1}\sharp B)][I\nabla_t (A^{-1}\sharp B)]^{-1} \le L_{s,t}  I.
  $$
 Note that $I\nabla_s (A^{-1}\sharp B)$ and $[I\nabla_t (A^{-1}\sharp B)]^{-1}$ are commuting, 
  and \eqref{BW geodesic form 2} implies that
  \begin{eqnarray*}
A\diamond_s B = [I\nabla_s (A^{-1}\sharp B)][I\nabla_t (A^{-1}\sharp B)]^{-1}(A\diamond_t B) [I\nabla_s (A^{-1}\sharp B)][I\nabla_t (A^{-1}\sharp B)]^{-1}.     
\end{eqnarray*}
Applying Theorem \ref{thm:Hermitian-near}, we get
\begin{equation}
   l_{s,t}^2 (A\diamond_t B) \preceq A\diamond_s B\preceq L_{s,t}^2 (A\diamond_t B).
\end{equation}
Obviously, $A\preceq B$ implies that $\mu_n\ge 1$, and $A\succeq B$ implies that $\mu_1\le 1$.
Using the expression \eqref{f_st}, it is straightforward to verify that 
$$
0< l_{s,t}=f_{s,t}(\mu_n)\le \cdots\le f_{s,t}(\mu_1)=L_{s,t}
$$
for the three cases given in Theorem \ref{thm:near_order_range}, so that   \eqref{ws:s-t-range} holds. 

\item 
For $t,s\in\R$, $(A^{-1}\sharp B)^t$ and $[I\nabla_s (A^{-1}\sharp B)]$ commute. Moreover, 
when $\mu>1$ and $s\in (\frac{1}{1-\mu},\infty)$, or $0<\mu<1$ and $s\in (-\infty,\frac{1}{1-\mu})$,  we have
$\frac{1-s+s\mu}{\mu^t}>0$. Now
$$
A\diamond_s B = [I\nabla_s (A^{-1}\sharp B)](A^{-1}\sharp B)^{-t}(A\natural_t B)[I\nabla_s (A^{-1}\sharp B)](A^{-1}\sharp B)^{-t} 
$$
in which $m_{s,t} I\le  [I\nabla_s (A^{-1}\sharp B)](A^{-1}\sharp B)^{-t}\le M_{s,t} I.$
Applying Theorem \ref{thm:Hermitian-near}, we get \eqref{sg-ws:s-t-range}.  \qedhere  
\end{enumerate}
\end{proof}

\begin{corollary}
    \label{thm:natural_near_order}
Let $A, B\in \P_n$. Let   $\mu_1\ge\cdots\ge\mu_n$ be the eigenvalues of $A^{-1}\sharp B$. The following statements hold: 
\begin{enumerate}
    \item  For $\mu\in (0,\mu_n]$, the parametric curves
$\{(\mu ^2 A)\natural_t B \mid t\in\R \}$ and $\{(\mu ^2 A)\diamond_t B \mid t\ge 0\}$ are monotonically increasing with respect to the near order. 
    \item  For $\mu\in [\mu_1,\infty )$, the parametric curves
$\{(\mu^2 A)\natural_t B\mid t\in\R \}$  and $\{(\mu ^2 A)\diamond_t B \mid t\le 1\}$ are monotonically decreasing with respect to the near order. 
\end{enumerate}
\end{corollary}

\begin{proof}
Let  $t<s$.  \eqref{natural:s-t-range} shows that for  $\mu\in (0,\mu_n]$, 
$$
\mu ^{2s-2t} (A\natural_t B)\le \mu_n^{2s-2t} (A\natural_t B)\preceq A\natural_s B,
$$
so that
\begin{eqnarray*}
(\mu ^2 A)\natural_t B= \mu ^{2-2t} (A\natural_t B) \preceq  \mu ^{2-2s} (A\natural_s B)=(\mu ^2 A)\natural_s B.
\end{eqnarray*}
Hence $\{(\mu ^2 A)\natural_t B \mid t\in\R \}$  is monotonically increasing   with respect to  the near order. By Theorem \ref{thm: A preceq B}, 
$\{(\mu ^2 A)\diamond_t B \mid t\ge 0\}$   is monotonically increasing    with respect to  the near order. 
The proof of second statement is analogous. 
\end{proof}

By Theorem~\ref{thm:wm} (2), when $t\in (0,1)$, $(A\diamond_t B)^{-1} = A^{-1}\diamond_t B^{-1}$ if and only if $A=B$. 
%We consider whether there exists a relation between them for any $A, B\in \P_n$. 
In general, they satisfy the following near order relation.

\begin{theorem}\label{thm:wm_inverse}
    Let $A, B\in \P_n$. Then for $t\in (0,1)$,
    \begin{eqnarray}
          (A^{-1}\diamond_t B^{-1})^{-1}  &\preceq &  A\natural_t B,
\\        (A^{-1}\diamond_t B^{-1})^{-1} &\preceq & A\diamond_t B,
    \end{eqnarray}
and either  of the equalities holds if and only if $A= B$.  
\end{theorem}

\begin{proof}
According to \eqref{BW geodesic form 2}, 
\begin{eqnarray*}
    (A^{-1}\diamond_t B^{-1})^{-1} 
    &=& [I\nabla_t (A\sharp B^{-1})]^{-1}A [I\nabla_t (A\sharp B^{-1})]^{-1}
    \\
    &=& [I\nabla_t (A^{-1}\sharp B)^{-1}]^{-1}A [I\nabla_t (A^{-1}\sharp B)^{-1}]^{-1}.
\end{eqnarray*}
The following L\"owner order relations exist for $t\in (0,1)$: 
$$
[I\nabla_t (A^{-1}\sharp B)^{-1}]^{-1} \le (A^{-1}\sharp B)^t \le  I\nabla_t (A^{-1}\sharp B) .
$$
By Corollary \ref{thm:near}, 
we get  $(A^{-1}\diamond_t B^{-1})^{-1} \preceq  A\natural_t B$ and
 $(A^{-1}\diamond_t B^{-1})^{-1} \preceq A\diamond_t B$.  Moreover,  either  of the equalities holds if and only if $A^{-1}\sharp B=I$, that is, $A=B$. 
\end{proof}

\begin{remark}
Theorem~\ref{thm:wm} (6) shows that for $A, B\in \P_n$ and $t\in [0,1]$,
    \[(A^{-1}\diamond_t B^{-1})^{-1} \ge (A^{-1}\nabla_t B^{-1})^{-1}.\]
By  Theorems~\ref{thm:near-main} and~\ref{thm:wm_inverse},  the Wasserstein mean, the spectral geometric mean, the harmonic mean, and the arithmetic mean have the following relations:
\begin{equation}\label{means-order-seq-1}
         (A^{-1}\nabla_t B^{-1})^{-1}\le (A^{-1}\diamond_t B^{-1})^{-1}
        \preceq  (A^{-1}\natural_t B^{-1})^{-1}=A\natural_t B 
        \preceq   A\diamond_t B  \le  A \nabla_t B.
\end{equation}
  Equivalently,
\begin{equation}\label{means-order-seq-2}
        (A \nabla_t B)^{-1} \le  (A \diamond_t B)^{-1} \preceq  (A\natural_t B)^{-1}=A^{-1}\natural_t B^{-1} \preceq  A^{-1}\diamond_t B^{-1} \le A^{-1}\nabla_t B^{-1}. 
\end{equation}
By \cite[Theorem 2]{DF23}, we conclude that any two means in the sequence \eqref{means-order-seq-1} or \eqref{means-order-seq-2} have
at least the near order relation. 
\end{remark}

By  \cite[Proposition 4, Remark 5]{DF23}, or by an analogous result of \cite[Theorem 3.6]{HK22}, 
if $A, B\in \P_n$  and $A\preceq B$, then
 $A^p\preceq B^p$ for $p\ge 1$;   if $A\le B$ or $\log A\le \log B$, then $A^p\preceq B^p$ for $p\ge 0$. 
Note that $A\preceq B$ if and only if $B^{-1}\preceq A^{-1}$.  %(see \cite[Proposition 7]{DF23}). 
We summarize  the results as follows and skip their proofs. 
%Theorem~\ref{thm:new_relation}. 

\begin{theorem}\label{thm:new_relation}
Let $A, B\in \P_n$. Suppose that one of the following holds:
\begin{enumerate}
    \item $A\preceq B$ and $p\ge 1$, or 
    \item  $A\le B$ and $p\ge 0$, or 
    \item $\log A\le \log B$ and $p\ge 0$.
\end{enumerate}
Then 
 $A^p\preceq B^p$ and $B^{-p}\preceq A^{-p}$. Moreover,
 \begin{enumerate}
     \item the parametric curves $\{A^p\diamond_t B^p\mid t\ge 0\}$ and $\{B^{-p}\diamond_t A^{-p}\mid t\ge 0\}$  are geodesics  monotonically increasing with respect to the near order;
      \item the parametric curve  $\{A^p\natural_t B^p\mid t\in\R\}$ %% and $\{B^{-p}\natural_t A^{-p}\mid t\in\R\}$  are 
      is monotonically increasing with respect to the near order.
 \end{enumerate}
\end{theorem}

\begin{remark}
    The statement ``$A\preceq B$ and $p\ge 1$ imply that  $A^p\preceq B^p$'' 
    %can be proved by Ando-Hiai inequality \cite{AH94}: if  $A,B\in\P_n$ and $t\in (0,1)$, then $(A^{-1}\sharp_t B)^r\prec_{\log} A^{-r}\sharp_t B^{r}$ for $r\in [0,1]$ and     $(A^{-1}\sharp_t B)^r\succ_{\log} A^{-r}\sharp_t B^{r}$ for $r\in [1,\infty)$. 
    can be viewed as a special case of   \cite[Theorem 6]{DF23}, which can  be applied to obtain other  monotonic curves with respect to $\preceq $.
\end{remark}

\iffalse
\begin{remark}
    By \cite{Fu87}, if $A\le B$ in $\P_n$, then
\begin{equation}
(A^{r/2} B^s A^{r/2})^{\alpha} \ge A^{(s+r)\alpha}
\end{equation}
for all $r, s\ge 0$, $\alpha\in (0,1)$, and $\alpha s -(1-\alpha) r\le 1$. 
Let $\alpha=1/2$ then we have: if $A\le B$, then
\begin{equation}
    A^{-r}\sharp  B^{s} \ge A^{(s-r)/2}
\end{equation}
for all $r \ge 0$ and $s\in [0,r+2]$. 
\end{remark} 
 
Theorem~\ref{thm:new_relation} can give the results of the Wasserstein power mean and the spectral geometric power mean relations with respect to the near order.
% \blue{theorem 4.8 and remark 4.9 below have been included in Theorem \ref{thm:new_relation}. }
\begin{corollary}
 Let $A, B\in \P_n$ and $p\geq 1$. If $A\diamond_t B\preceq  B$ for some $t\in [0,1)$ or $A\preceq A\diamond_t B$ for some $t\in (0,1]$, then both of the following hold:
 \begin{enumerate}
     \item $A^p \preceq A^p \diamond_t B^p \preceq  B^p$ for all $t\in (0,1)$;
     \item $A^p \preceq A^p \natural_t B^p \preceq  B^p$ for all $t\in (0,1)$.
 \end{enumerate}
\end{corollary}

\begin{proof}
    When $A\diamond_t B\preceq  B$ for some $t\in [0,1)$ or $A\preceq A\diamond_t B$ for some $t\in (0,1]$, according to Theorem~\ref{thm: A preceq B}, it is straightforward to have $A\preceq B$, that is, $A^{-1}\sharp  B\geq I$.
    For all $p\ge 1$, by the Ando-Hiai inequality, we have $A^{-p}\sharp B^{p} \le I$, which is $A^p \preceq B^p$.
    By Theorem~\ref{thm: A preceq B} again, $A^p \diamond_t B^p \preceq  B^p$ for all $t\in [0,1)$.
\end{proof}
\fi

Because the L\"owner order is stronger than the near order and the near order is stronger than the eigenvalue entrywise order, it is straightforward to have the following corollary.

\begin{corollary}\label{thm:wm_power}
    Let $A, B\in \P_n$ and $p\geq 1$. If $A\diamond_t B\preceq A\diamond_s B$ for some $0\le t<s\le 1$, or 
     $A\natural_t B\preceq A\natural_s B$ for some $ t<s $, then the  matrices on the curves
    $\{A^p \diamond_t B^p\mid t\in [0,\infty)\}$ and  $\{A^p\natural_t B^p\mid t\in\R\}$ have entrywise monotonically increasing eigenvalues. In particular, 
    \begin{enumerate}
        \item $A^p \le_{\la} A^p \diamond_t B^p \le_{\la} B^p$ for all $t\in (0,1)$,
        \item $A^p \le_{\la}A^p\natural_t B^p \le_{\la} B^p$ for all $t\in (0,1)$.
    \end{enumerate}
\end{corollary}

% \begin{proof}
% Since $A\diamond_t B\le B$ for some $t\in [0,1)$, by Lemma~\ref{lem:HK22}, $B^{-p}\sharp  A^p \le I$ for all $p\ge 1$.
% Then by Lemma~\ref{lem:power}, we obtain
% \[[I\nabla_{1-t}(B^{-p}\sharp  A^p)]^2 \le [I\nabla_{1-t} I]^2 = I.\]
% By~\eqref{eqn:Wasserstein}, we obtain
% \[ A^p\diamond_t B^p=B^p\diamond_{1-t} A^p=_{\la}B^{p/2}[I\nabla_{1-t}(B^{-p}\sharp  A^p)]^2 B^{p/2} \le B^{p/2}I B^{p/2} = B^p.\]
% The proof is completed.
% \end{proof}

% Similarly, by Lemma~\ref{lem:AAB}, we can obtain the following theorem. The proof is identical to the proof of Theorem~\ref{thm:wm_power} so we skip the proof here.

% \begin{corollary}\label{thm:wm_power 2}
        % Let $A, B\in \P_n$. If $A\le A\diamond_t B$ for some $t\in (0,1]$, then $A^p \le_{\la}A^p\diamond_t B^p $ for all $p\geq 1$.    
% \end{corollary}

Corollary \ref{thm:wm_power} can be applied to the following two cases.

\begin{corollary}
       Let $A, B\in \P_n$.  If $A \ge I$ and $A\diamond_t B \leq I$ for some $t\in(0, 1)$, then \[ A^p \diamond_s B^p \leq \left(\frac{1}{1-t+t\lambda_n(A^{-1}\sharp  B)}\right )^{2p} I \] 
       for all $p\ge 1$ and all $s\in [0,1]$.     
\end{corollary}

\begin{proof}
    % Note that
    % $$A\diamond_t B=[I\nabla_t(A^{-1}\sharp  B)]A[I\nabla_t(A^{-1}\sharp  B)]. $$
    Since $A\diamond_t B \leq I$, according to~\eqref{eqn:Wasserstein}, 
    \[A\leq [I\nabla_t(A^{-1}\sharp  B)]^{-2} \le (1-t+t\lambda_n(A^{-1}\sharp  B) )^{-2} I.\]
    Since $A\diamond_t B\le I\le A$, by Corollary~\ref{thm:wm_power}, for all $p\ge 1$ and all $s\in [0,1]$,
    \[\la_1(A^p \diamond_s B^p) \le \la_1(A^p) \le \left(\frac{1}{1-t+t\lambda_n(A^{-1}\sharp  B)}\right )^{2p},\]
    which means that $A^p \diamond_s B^p \leq \left(\frac{1}{1-t+t\lambda_n(A^{-1}\sharp  B)}\right )^{2p} I$.
\end{proof}

A similar argument on $A^p \natural_s B^p$ leads to the following result. We skip its proof here.

\begin{corollary}
         Let $A, B\in \P_n$. If $A \ge I$ and $A\natural_t B \leq I$ for some $t\in(0, 1)$, then \[ A^p \natural_s B^p \leq \left(\frac{1}{\lambda_n(A^{-1}\sharp  B)}\right )^{2pt} I\] 
       for all $p\ge 1$ and all $s\in [0,1]$. 
\end{corollary}

Using $A\diamond_t B=B\diamond_{1-t} A$ and $A\natural_t B=B\natural_{1-t} A$, Theorem \ref{thm:new_relation} can be changed to an equivalent theorem for curves monotonically decreasing with respect to the near order.

\begin{theorem}\label{thm:new_relatio_2}
Let $A, B\in \P_n$. Suppose that one of the following holds:
\begin{enumerate}
    \item $A\succeq B$ and $p\ge 1$, or 
    \item  $A\ge B$ and $p\ge 0$, or 
    \item $\log A\ge \log B$ and $p\ge 0$.
\end{enumerate}
Then 
 $A^p\succeq B^p$ and $B^{-p}\succeq A^{-p}$. Moreover,
 \begin{enumerate}
     \item the parametric curves $\{A^p\diamond_t B^p\mid t\le 1\}$ and $\{B^{-p}\diamond_t A^{-p}\mid t\le 1\}$  are geodesics  monotonically decreasing with respect to the near order;
      \item the parametric curve  $\{A^p\natural_t B^p\mid t\in\R\}$ %% and $\{B^{-p}\natural_t A^{-p}\mid t\in\R\}$  are 
      is monotonically decreasing with respect to the near order.
 \end{enumerate}
\end{theorem}

% \red{Include the comparison to the curve $A\natural_t B$ for $t\ge 0$ using $$A^{-1} \sharp  (A\natural_t B) \leq A^{-1} \sharp (A\diamond_t B).$$
% Corollary: When $A\preceq A\natural_t B$ for certain $t\in (0,1)$, we  have $A\preceq B$.
% Conversely, if $A\preceq B$, then $A\natural_t B\preceq B$ for all $t\in (0,1)$.  Theorem \ref{thm: A preceq B} can be applied to these cases.
% }

% We may extend the definition of $A\natural_t B$ to all $t\ge 0$ and show that for $t\ge 0$: 
% $$A^{-1}\sharp (A\natural_t B)=(A^{-1}\sharp  B)^t\le I\nabla_t (A^{-1}\sharp  B)= A^{-1}\sharp  (A\diamond_t B).$$
% Also $A\preceq B$ if and only if $A\preceq A\natural_t B$ for certain $t\ge 0$ or $A\natural_t B\preceq B$ for certain $t\in (0,1).$

% Note that $X=A^{-1}\sharp  C$ is the unique solution in $P_n$ to $XAX=C$. We get $A\diamond_t B=X_tAX_t$ for
% a curve $\{X_t\mid t\ge 0\}\subseteq  P_n$ where $t\le s$ implies $X_t\le X_s$.

\section{The L\"owner orders on two means}\label{Sec:properties}

% \begin{theorem}\cite{ZH23} \label{thm:zheng}
%     For $A, B\in \P_n$ and $t\in [0,1]$, 
%     \[A\diamond_t B =B\diamond_{1-t} A = |(1-t)A^{1/2} +t U^*B^{1/2} |^2,\]
%     in which $U$ is certain unitary matrix. Moreover, when 
%     % $(1-t)A^{1/2} +t |B^{1/2}A^{1/2}| \geq 0$ 
%     $t\in [0,1]$,
%     \[|(A\diamond_t B)^{1/2}A^{1/2}| = (1-t)A+t |B^{1/2}A^{1/2}|.\]
% \end{theorem} 

We explore the L\"owner order properties of Wasserstein mean and spectral geometric mean in this section. These  properties
further display the interesting similarities  between the two means. 
 Some of them can be extended to the curves induced by the two means. 
 
%Let us start with some identities.  
\eqref{BW geodesic form 2}, \eqref{spectral geometric form}, and Theorem \ref{thm:gm}(9) 
imply that for $A,B\in\P_n$ and $t\in [0,1]$: 
\begin{eqnarray}\label{inv-sharp-diamond}
    A^{-1}\sharp (A\diamond_t B)=B\sharp (B^{-1}\diamond_t A^{-1}) &=& I\nabla_t (A^{-1}\sharp B),
\\\label{inv-sharp-natural}
A^{-1}\sharp (A\natural_t B)=B\sharp (B^{-1}\natural_t A^{-1}) &=&  (A^{-1}\sharp B)^t.
\end{eqnarray}
They can be used to derive  identities like:
\begin{eqnarray}
    [A\sharp (A\diamond_t B)^{-1} ] \nabla_t [ B\sharp (A\diamond_t B)^{-1} ] &=& I,
 \\  \ [A\sharp (A\natural_t B)^{-1} ]^{1-t} [ B\sharp (A\natural_t B)^{-1} ]^{t} &=& I. 
\end{eqnarray}
\eqref{inv-sharp-diamond} and \eqref{inv-sharp-natural} also imply the following  L\"owner order relations.

\begin{theorem}\label{thm:order_sp_was}
    For $A, B\in \P_n$ and $t\in[0,1]$,  
\begin{eqnarray}
    \label{eqn:order_sp_was}
    A^{-1} \sharp  (A\natural_t B) &\leq & A^{-1} \sharp (A\diamond_t B),
    \\
     \label{eqn:order_sp_was 2}
    B^{-1} \sharp  (A\natural_t B) &\leq & B^{-1} \sharp (A\diamond_t B),
\end{eqnarray}
\end{theorem}

\eqref{eqn:order_sp_was} and \eqref{eqn:order_sp_was 2} can be rephrased as follows: 
\[|(A\natural_t B)^{1/2}A^{1/2}|\leq |(A\diamond_t B)^{1/2}A^{1/2}|,\qquad
|(A\natural_t B)^{1/2}B^{1/2}|\leq |(A\diamond_t B)^{1/2}B^{1/2}|.\]

The above inequlities lead to  the following results.

\begin{theorem}\label{thm:was_geom}
   Let $A, B, C\in  \P_n$. 
\begin{enumerate}
    \item  If $A\diamond_t B\le A\diamond_t C$ or $A\natural_t B\le A\natural_t C$ for one $t\in(0, 1]$,  then $A^{-1}\sharp _s B \le A^{-1}\sharp _s C$ for all $s\in [0,1/2]$. 
    \item   If $A\diamond_t B=A\diamond_t C$  or $A\natural_t B =A\natural_t C$ for one $t\in(0, 1]$,  then $B=C$. 
\end{enumerate}
\end{theorem}

\begin{proof} We prove the statements related to $\diamond $, and the proofs for statements related to $\natural $ are analogous. 
\begin{enumerate}
    \item By \eqref{inv-sharp-diamond} and the assumption $A\diamond_t B\le A\diamond_t C$,
    \begin{eqnarray*}
A^{1/2}[I\nabla_t (A^{-1}\sharp  B)] A^{1/2}  &=&  (A^{1/2}(A\diamond_t B) A^{1/2} )^{1/2}
    \\
    &\le & ( A^{1/2}(A\diamond_t C) A^{1/2} )^{1/2}\\
    &=& A^{1/2}[I\nabla_t (A^{-1}\sharp  C)] A^{1/2} ,
    \end{eqnarray*}
    which gives $A^{-1}\sharp  B \le A^{-1}\sharp  C$. In other words, 
    $$(A^{1/2}BA^{1/2})^{1/2}\le (A^{1/2}CA^{1/2})^{1/2}.$$
    Therefore, for any $s\in [0,1/2],$ we have $2s\in [0,1]$ so that by L\"owner–Heinz inequality (see~\cite{Z02})
    $$(A^{1/2}BA^{1/2})^{s}\le (A^{1/2}CA^{1/2})^{s}.$$
    In other words, $A^{-1}\sharp _s B \le A^{-1}\sharp _s C$. 
    
    \item Similar to the above analysis, $A\diamond_t B=A\diamond_t C$ leads to $A^{-1}\sharp  B =A^{-1}\sharp  C$, which gives $B=C$.   \qedhere
\end{enumerate}
%Analogously, \eqref{eqn:sm} implies that $A^{1/2}(A^{-1}\sharp  B)^t A^{1/2}  = (A^{1/2}(A\natural_t B) A^{1/2} )^{1/2}$, which leads to proofs of the statements for $\natural$. 
\end{proof}

Theorem \ref{thm:was_geom} implies the following result, in which $A^{-1}\diamond_t A=(A^{-1/2}\nabla_t A^{1/2})^2$
for $A\in\P_n$. 

\begin{theorem}\label{equiv conds of A le B}
    Let $A, B\in \P_n$. Then $A\le B$ if one of the following holds for any $t\in (0,1)$:
\begin{enumerate}
    \item $A^{-1}\diamond_t A \le A^{-1}\diamond_t B$,
    \item $B^{-1}\diamond_t B\le A^{-1}\diamond_t B$,
    \item $B^{-1}\diamond_t A\le A^{-1}\diamond_t A$,
    \item $B^{-1}\diamond_t A\le B^{-1}\diamond_t B$.
\end{enumerate}
\end{theorem}
 
\begin{proof}
We prove the first statement and the others can be done similarly.

If $ A^{-1}\diamond_t A\le A^{-1}\diamond_t B $ for any $t\in (0,1)$, then by 
Theorem \ref{thm:was_geom}, 
$A=A \sharp  A\le A \sharp  B,$  so that
$I\le |B^{1/2}A^{-1/2}|$ and thus $A\le B$.
\end{proof}

The analogous results for the spectral geometric mean are given  as follows and we omit the proofs here. 
Note that $A^{-1}\natural_t A=A^{2t-1}$ for $A\in\P_n$. 

\begin{theorem}\label{equiv conds of A le Bspectral}
    Let $A, B\in \P_n$. Then $A\le B$ if one of the following holds for $t\in (0,1)$:
\begin{enumerate}
    \item $A^{-1}\natural_t A\le A^{-1}\natural_t B$,
    \item $B^{-1}\natural_t B\le A^{-1}\natural_t B$,
    \item $B^{-1}\natural_t A\le A^{-1}\natural_t A$,
    \item $B^{-1}\natural_t A\le B^{-1}\natural_t B$.
\end{enumerate}
\end{theorem}

The identity matrix $I$ commutes with any matrix $X\in\P_n$, so that for $t\in (0,1)$: 
\begin{equation}
   I\nabla_t X\ge  I\diamond_t X=(I\nabla_t X^{1/2})^2\ge X^t=I\natural_t X 
\end{equation}
and either equality holds if and only if $X=I.$

The following result is only for the Wasserstein mean. %has no   counterpart for the spectral geometric mean yet. 

\begin{theorem}
    If $A, B\in \overline{\P_n}$, then for $t\in (0,1)$,
    \[|(A\diamond_t B)^{1/2}A^{1/2}| \ge  A\diamond_t |B^{1/2}A^{1/2}|.\]
\end{theorem}
\begin{proof}
By~\eqref{BW geodesic form},  we have
\begin{equation*}
    |(A\diamond_t B)^{1/2}A^{1/2}|  =  (1-t)A+t|B^{1/2}A^{1/2}|= A\nabla_t |B^{1/2}A^{1/2}| \geq A\diamond_t |B^{1/2}A^{1/2}|.\qedhere
\end{equation*}
\end{proof}

\bibliographystyle{abbrv}
\bibliography{refs} % used citations in bib.bib have been put in refs.bib 

\end{document}